\def\ngay{April 3, 2008}
\documentclass[12pt,a4]{amsart}
\usepackage{amsmath,amssymb,amsfonts,amscd}
\usepackage[all]{xy}
\usepackage{ifthen}
\usepackage{hyperref}

\begin{document}

\theoremstyle{plain}
\newtheorem{thm}[subsection]{Theorem}
\newtheorem{lem}[subsection]{Lemma}
\newtheorem{cor}[subsection]{Corollary}
\newtheorem{prop}[subsection]{Proposition}
\newtheorem{remark}[subsection]{Remark}
\newtheorem{defn}[subsection]{Definition}
\newtheorem{ex}[subsection]{Example}
\newtheorem{reduct}[subsection]{Reduction}
\newtheorem{conj}[subsection]{Conjecture}

\numberwithin{equation}{section}

\newcommand{\otimesover}[3][]{\ifthenelse{\equal{#1}{}}%
{\xymatrix@1@=0pt@M=0pt{_{#2}&\otimes\hspace{1pt} &_{#3}}}%
{\xymatrix@1@=0pt@M=0pt{_{#2}&\otimes\hspace{1pt} &_{#3}\\& ^{#1}&}}}

\newcommand{\timesover}[3][]{\ifthenelse{\equal{#1}{}}%
{\xymatrix@1@=0pt@M=0pt{_{#2}&\times\hspace{1pt} &_{#3}}}%
{\xymatrix@1@=0pt@M=0pt{_{#2}&\times\hspace{1pt} &_{#3}\\& ^{#1}&}}}

\newcommand{\sqtimesover}[3][]{\ifthenelse{\equal{#1}{}}%
{\xymatrix@1@=0pt@M=0pt{_{#2}&\Box\hspace{1pt} &_{#3}}}%
{\xymatrix@1@=0pt@M=0pt{_{#2}&\Box\hspace{1pt} &_{#3}\\& ^{#1}&}}}

\newcommand{\eq}[2]{\begin{equation}\label{#1}#2 \end{equation}}
\newcommand{\ml}[2]{\begin{multline}\label{#1}#2 \end{multline}}

\def\printlabel{0}
\newcommand{\ga}[2]{\ifthenelse{\printlabel=1} 
{\marginpar{#1}\begin{gather}\label{#1}#2\end{gather}}
{\begin{gather}\label{#1}#2 \end{gather}}
}
\let\supportlabel\label
\marginparwidth10ex
 \def\label#1{\ifthenelse{\printlabel=1}
  {\protect \supportlabel{#1}\text{\underline{#1}}}
  {\supportlabel{#1}}
}

\newcommand{\surj}{\twoheadrightarrow}
\newcommand{\inj}{\hookrightarrow}
\newcommand{\Ker}{{\rm Ker\hspace{.1ex}}}
\newcommand{\rank}{{\rm rank}}
\newcommand{\Hom}{{\rm Hom}}
\newcommand{\im}{{\rm im}}
\newcommand{\Spec}{{\rm Spec\hspace{.1ex}}}
\newcommand{\Tr}{{\rm Tr}}
\newcommand{\Gal}{{\rm Gal}}
\newcommand{\sA}{{\mathcal A}}
\newcommand{\sB}{{\mathcal B}}
\newcommand{\sC}{{\mathcal C}}
\newcommand{\sD}{{\mathcal D}}
\newcommand{\sE}{{\mathcal E}}
\newcommand{\sF}{{\mathcal F}}
\newcommand{\sG}{{\mathcal G}}
\newcommand{\sH}{{\mathcal H}}
\newcommand{\sI}{{\mathcal I}}
\newcommand{\sJ}{{\mathcal J}}
\newcommand{\sK}{{\mathcal K}}
\newcommand{\sL}{{\mathcal L}}
\newcommand{\sM}{{\mathcal M}}
\newcommand{\sN}{{\mathcal N}}
\newcommand{\sO}{{\mathcal O}}
\newcommand{\sP}{{\mathcal P}}
\newcommand{\sQ}{{\mathcal Q}}
\newcommand{\sR}{{\mathcal R}}
\newcommand{\sS}{{\mathcal S}}
\newcommand{\sT}{{\mathcal T}}
\newcommand{\sU}{{\mathcal U}}
\newcommand{\sV}{{\mathcal V}}
\newcommand{\sW}{{\mathcal W}}
\newcommand{\sX}{{\mathcal X}}
\newcommand{\sY}{{\mathcal Y}}
\newcommand{\sZ}{{\mathcal Z}}
\newcommand{\A}{{\Bbb A}}
\newcommand{\B}{{\Bbb B}}
\newcommand{\C}{{\Bbb C}}
\newcommand{\D}{{\Bbb D}}
\newcommand{\E}{{\Bbb E}}
\newcommand{\F}{{\Bbb F}}
\newcommand{\G}{{\Bbb G}}
\renewcommand{\H}{{\Bbb H}}
\newcommand{\I}{{\Bbb I}}
\newcommand{\J}{{\Bbb J}}
\newcommand{\bl}{{\mathbb L}}
\newcommand{\M}{{\Bbb M}}
\newcommand{\N}{{\Bbb N}}
\renewcommand{\P}{{\Bbb P}}
\newcommand{\Q}{{\Bbb Q}}
\newcommand{\R}{{\mathbb R}}
\newcommand{\T}{{\Bbb T}}
\newcommand{\U}{{\Bbb U}}
\newcommand{\V}{{\Bbb V}}
\newcommand{\W}{{\Bbb W}}
\newcommand{\X}{{\Bbb X}}
\newcommand{\Y}{{\Bbb Y}}
\newcommand{\Z}{{\Bbb Z}}

\newcommand{\fg}{{\mathfrak g}}
\newcommand{\fh}{{\mathfrak h}}
\newcommand{\fp}{{\mathfrak p}}
\newcommand{\fq}{{\mathfrak q}}
\newcommand{\ff}{{\mathfrak f}}

\newcommand{\Mod}{\text{\sf Mod}}
\newcommand{\vect}{\text{\sf vect}}
\newcommand{\Vect}{\text{\sf Vect}}
\newcommand{\Rep}{\text{\sf Rep}}
\newcommand{\id}{{\rm id\hspace{.1ex}}}
\newcommand{\Ind}{\text{\rm Ind-}}
\def\Sp{\Spec}
\newcommand{\res}{{\text{\sf res}\hspace{.1ex} }}
\newcommand{\ind}{{\text{\sf ind}\hspace{.1ex}}}

\title[Construction of quotient tensor category]{A construction of a quotient tensor category}

\author{Ph\`ung  H\^o Hai}
\address{
Universit\"at Duisburg-Essen, Essen, Germany  and
Hanoi Institute of Mathematics, Hanoi, Vietnam}
\email{hai.phung@uni-duisburg-essen.de \& phung@math.ac.vn}
\date\ngay

\begin{abstract} 
For a rigid tensor abelian category $\sT$ over a field $k$ we introduce a notion
of a normal quotient $\fq:\sT\longrightarrow\sQ$. 
In case $\sT$ is a Tannaka category, our notion is equivalent to Milne's
notion of a normal quotient \cite{milne2005}. More precisely,
if $\sT$ is  the category of finite dimensional representations of a groupoid
scheme $G$ over $k$, then $\sQ$ is equivalent to the representation category of
a normal subgroupoid scheme of $G$. We describe such a quotient
in terms of the subcategory $\sS$ of $\sT$ consisting of objects which become
trivial in $\sQ$. We show that, under some condition on $\sS$, $\sQ$ is uniquely
determined by $\sS$.
If $\sS$ is a finite tensor category, we show that the quotient of $\sT$ by $\sS$ exists.
In particular we show the existence of the base change of $\sT$ with respect to finite field extensions.  As an application, we obtain a condition for the exactness of sequences of groupoid schemes in terms of the representation categories.
\end{abstract}

\maketitle
\section*{Introduction}
Let $\sT$ be a Tannaka category over $k$ with fiber functor 
$\omega$ to $\vect_K$, where $K\supset k$ a field extension and let
$\sS$ be a full tensor abelian subcategory of $\sT$ which is  closed
under taking subquotients. The natural inclusion $\sS\longrightarrow \sT$ induces
a surjective homomorphism of $k$-groupoid schemes 
(hereafter called groupoids for short) acting transitively upon $\Sp (K)$
\ga{eq2}{G(\sT)\longrightarrow G(\sS).}
The kernel of this homomorphism is
a $K$-group scheme $L$. Denote by 
$\sQ$ the category of finite dimensional representations of $L$ 
and by $\fq:\sT\longrightarrow \sQ$ the restriction functor. 
For the sequence of functors
\ga{eq3}{\sS\longrightarrow \sT\stackrel{\fq}{\longrightarrow} \sQ} 
it is shown in \cite{EsnHai05} that
\begin{itemize}
\item[] each object of $\sQ$ is isomorphic (in $\sQ$) to
a subobject of the image under $\fq$ of an object of $\sT$.\end{itemize}

The problem we want to address in this work is to give an abstract description of 
the quotient category $\sQ$. This question was posed by P. Deligne in connection with our
description of the representation category of $L$ given in \cite{EsnHai05}.

A related question was studied by J. Milne \cite{milne2005}, where the
normal quotients $(\sQ,q)$ were classified in terms of $\sS$. We notice that
Milne introduces the notion of normal quotients only for Tannaka categories and he
considers only quotients by a neutral Tannaka subcategory.

While studying this problem we realized that the Tannaka assumption on $\sT$
can be replaced by a weaker assumption. In the sequel we will assume that
$\sT$ is a rigid tensor abelian 
category over $k$ (see \ref{sect-def01} for precise definition), which does not necessarily
have a fiber functor. 
Let us start with our definition of
a normal quotient category of such a category $\sT$.

A normal $K$-\emph{quotient} of $\sT$ is by definition a pair $(\sQ,\fq:\sT\longrightarrow\sQ)$
consisting of a rigid tensor  abelian category $\sQ$ over $K$ and a 
 $k$-linear exact tensor functor $\fq$ (the $k$-linear structure over
$\sQ$ is induced from the inclusion $k\subset K$), such that:
\begin{itemize}\item[(i)] for an object $X$  of $\sT$
the largest trivial subobject of $\fq(X)$  is isomorphic to the image
under $\fq$ of a subobject of $X$, (a trivial object is a direct sum of copies of the
unit object);
\item[(ii)] each object of $\sQ$ is isomorphic to a subobject of $\fq(X)$ with $X\in\sT$.
\end{itemize}
We notice here that our condition ``subobject'' in (ii) above is
stronger than the usual condition ``subquotient'' used in \cite{DeMil}.
This condition was first mentioned in \cite{EsnHai05}. It reflects the difference
between subgroups and normal subgroups.

Given a $K$-quotient $(\sQ,\fq)$ of $\sT$, let $\sS$ denote the full
subcategory of $\sT$ consisting of those objects which become trivial in $\sQ$. Then
$\sS$ is a tensor abelian subcategory and is closed under taking subquotients.
We shall call $\sS$ the invariant subcategory with respect to
the quotient $(\sQ,\fq)$. 
For each object $X$ of $\sT$, let $X_\sS$ denote the maximal $\sS$-subobject of $X$,
i.e. the maximal subobject of $X$ which is isomorphic to an object of $\sS$.
The restriction of $\fq$ to $\sS$ can be considered as a fiber functor to $\vect_K$ (Lemma \ref{lem-def03}).
We shall denote this fiber functor by $\omega$.

\subsection{}We first consider the case when $K=k$. Thus $\omega$ is a neutral fiber functor for $\sS$. 
Our main assumption is that $\sT$ is flat over $\sS$ (see Definition \ref{def-q08}).
With this assumption, our first result (see \ref{sect-q10})
states that
\begin{itemize}\item[(A)]\em  the category $\sQ$ 
is equivalent to the category, whose objects are triples
$(X,Y,f\in\omega(X^\vee\otimes Y)_\sS)$, which has to be understood as the image
of $X$ in $Y$ by $f$, and morphisms are 
appropriately defined. Consequently $(\sQ,\fq)$ is
 uniquely determined by $(\sS,\omega)$, up to a tensor equivalence.
\end{itemize}
Thus, in this case we can say that $\sQ$ is a quotient of $\sT$ by $\sS$.
Further we show that (Corollary \ref{cor-q12})
\begin{itemize}\item[(B)] \em if $\sT$ is a Tannaka category then the quotient of $\sT$ by any
neutral Tannaka category, which is closed under taking subquotients, exists.
\end{itemize}
This result shows in particular that our notion of normal quotients,
when restricted to Tannaka categories, coincides with Milne's definition. 
Unfortunately we cannot establish the existence of the quotient in the general case.
However, when $\sS$ is an \'etale finite category, we show in Theorem \ref{thm-q18} that
\begin{itemize}\item[(C)] \em the quotient of $\sT$ by an \'etale finite neutral Tannaka
full subcategory, closed under taking subquotients, exists.\end{itemize}

\subsection{}
To treat the general case ($K\neq k)$ we first need the notion of base change for $\sT$.
By definition, the base change $\sT_{(K)}$ of $\sT$ is a normal $K$-quotient
of $\sT$, the invariant category of which consists of trivial objects.
Under the assumption that $\sT$ is flat over $K$ (see \ref{b03}) we describe the
base change $\sT_{(K)}$ of $\sT$.
Again the existence of base change for Tannaka category is known, see \ref{cor-b06}, but it is not
known for arbitrary rigid tensor abelian categories. However, we are able to
establish the existence of base change with respect to finite separable field extensions, see \ref{sect-b08}.

\subsection{}Now, assuming that $\sT$ is flat over $K$, that the base change $\sT_{(K)}$ exists
and that $\sT$ is flat over $\sS$ (see  \ref{def-g08}) we show a similar result to (A) above
(Theorem \ref{thm-g09}). The existence of the quotient by an \'etale finite tensor full subcategory
(closed under taking subquotients) is also established (cf. \ref{sect-g04}).
\subsection{}
Our description of quotient categories, when applied to the representation categories of
groupoid schemes, yields a criterion for the exactness of sequences of groupoid schemes
(Corollary \ref{cor-g11}).

As mentioned above, we cannot show the existence of the quotient of $\sT$ by
a Tannaka full subcategory, closed under taking subquotients. However
the existence for the case $\sS$ is \'etale finite gives evidence for believing that
the quotient might exist in a more general case. More precise questions are given in
\ref{rm-g12}.

\bigskip

\begin{center}{\bf Acknowlegment}\\
\end{center}
The author thanks Prof. H. Esnault for  communicating with him
the question and suggestion of Prof. P. Deligne and for many stimulating 
and enlightening discussions. This
work results from these discussions. The author thanks
Prof. P.~Deligne for many useful comments which greatly
improves presentation of the work.
The work was partially supported by a
Heisenberg-Stipendium of the DFG,
the Leibniz-Preis of Professors H.~Esnault
and E.~Viehweg and the NFBS of Vietnam.

\section{Preliminaries}\label{sect-pre}
We consider the category of affine group schemes (not necessarily of finite type)
 over a field $k$, which
 we shall call groups for short. Let $L$ denote the kernel of
a homomorphism $f:G\longrightarrow A$, which is then
 a normal subgroup of $G$. We collect here some known information on the relationship
 between these groups.

 \subsection{}\label{sect-pre01}
 We shall use the notation $\sO(G)$ to denote the function algebra over $G$, which
 is a Hopf $k$-algebra. The reader is referred to \cite{waterhouse} for details on
the structure of $\sO(G)$.  Let $m,u,\epsilon,\Delta,\iota$ denote the product,
unit map, coproduct, counit map and the antipode of $\sO(G)$, respectively.
We shall adopt Sweedler's notation for the coproduct: 
$$\Delta(h)=\sum_{(h)}h_{(1)}\otimes h_{(2)}.$$

\subsection{}\label{sect-pre02}
The category $\Rep(G)$ of $k$-linear representation of 
$G$ is equivalent to the category $\sO(G)$-Comod of  $\sO(G)$ comodules,
which consists of pairs $(V,\rho_V:V\longrightarrow V\otimes \sO(G))$, where $V$ is
a vector space and $\rho_V$ is a $k$-linear map satisfying the following conditions:
\begin{eqnarray*}
(\rho_V\otimes\id)\rho_V=(\id\otimes\Delta)\rho_V:V\longrightarrow V\otimes V\otimes\sO_G,&&
(\id\otimes\epsilon)\rho_V=\id_V.
\end{eqnarray*}

The coproduct $\Delta:\sO(G)\longrightarrow \sO(G)\otimes\sO(G)$ defines a 
comodule of $\sO(G)$ on itself, this comodule corresponds to the right
 regular representation of $G$ in $\sO(G)$. The terminologies
$G$-equivariant and $\sO(G)$-colinear are equivalent.

A homomorphism $f:G\longrightarrow A$ is the same as a homomorphism of Hopf algebras
$f^*:\sO(A)\longrightarrow\sO(G)$. The fundamental theorem of algebraic group theory claims
that $\sO(G)$ is faithfully flat over its subalgebra $f^*(\sO(A))$, cf.  \cite{GaDem}.

\subsection{}\label{sect-pre03}Let $q:L\longrightarrow G$ be the kernel of $f$. 
 That is, $L$ is the fiber  of $f$  at the unit element of $A$, $e:\Spec k\longrightarrow A$.
Thus we have
\ga{pre01}{\sO(L)\cong \sO(G)\otimes_{\sO(A)}k.
}
 The homomorphism $q^*:\sO(G)\longrightarrow \sO(L)$ is just the projection
$\sO(G)\longrightarrow \sO(G)\otimes_{\sO(A)}k$ obtained by tensoring $\sO(G)$ with
the counit map $\epsilon=e^*:\sO(A)\longrightarrow k$.

 We shall assume from now on that $\sO(A)$ is a Hopf subalgebra of
$\sO(G)$ and $f^*$ is the inclusion.
 
\subsection{}\label{sect-pre04} The homomorphism $q:L\longrightarrow G$ induces 
a tensor functor $\res^q:\Rep(G)\longrightarrow \Rep(L)$,
which restricts a representation of $G$ in a vector space $V$
through $q$ to a representation of $L$. 
The functor $\res^q:\Rep(G)\longrightarrow \Rep(L)$ admits a right adjoint, which is the
induced representation functor $\ind_q:\Rep(L)\longrightarrow \Rep(G)$, that is we have
a functorial isomorphism
\ga{pre02}{\Hom_L(\res^q(V),U)\cong \Hom_G(V,\ind_q(U)),\quad V\in\Rep(G), U\in\Rep(L).
}
The functoriality yields a canonical $L$-linear map 
$$\varepsilon_U:\res^q\ind_q(U)\longrightarrow U.$$
The functor $\ind_q$ can be explicitly computed in terms of the invariant space
functor $(-)^L$. We prefer here the following Hopf algebraic description, which will
be exploited in the next sections.

For an $L$-representation $U$, 
denote by $U\Box_{\sO(L)}\sO(G)$ the equalizer of the following maps
\ga{pre03}{\xymatrix{
\ar[rrd]^{\rho_U\otimes\id}
U\otimes\sO(G)\ar[d]_{\id\otimes\Delta}& \\ 
U\otimes\sO(G)\otimes\sO(G)
\ar[rr]_{\id\otimes q^*\otimes\id}&&U\otimes\sO(L)\otimes\sO(G).
}}
and call it the cotensor product over $\sO(L)$
of $U$ with $\sO(G)$. Using Sweedler's notation we can describe $U\Box_{\sO(L)}\sO(G)$
as the set of elements $ \sum_iu_i\otimes g_i\in U\otimes \sO(G) $, such that
\ga{}{\notag 
\sum_{i,(u_i)} u_{i(0)}\otimes u_{i(1)}\otimes g_i=
\sum_{i,(g_i)}u_i\otimes q^*(g_{i(1)})\otimes g_{i(2)}.}

We have a functorial isomorphism
\ga{pre04}{\ind_q(U)\cong U\Box_{\sO(L)}\sO(G).}
For the cotensor product there is the following key isomorphism first considered
by Takeuchi \cite{takeuchi79}
\ga{pre05}{(U\Box_{\sO(L)}\sO(G))\otimesover[\sO(A)]{}{}\sO(G)\cong U\otimes \sO(G), \quad
\sum_i u_i\otimes g_i\otimes h\mapsto \sum_{i}u_i\otimes  g_ih,}
which together with the faithful flatness of $\sO(G)$ over $\sO(A)$ (cf. \ref{sect-pre02})
 show that the functor $\ind_q=(-)\Box_{\sO(L)}\sO(G)$
is faithfully exact. 
A direct consequence of \eqref{pre05} is the following equality
\ga{pre06}{k\Box_{\sO(L)}\sO(G)= \sO(A).}

For  a representation $V$ of $G$, we have the following $G$-equivariant isomorphism
\ga{pre07}{V\otimes \sO(G)\longrightarrow (V)\otimes\sO(G):=\sO(G)^{\oplus\dim_kV}, \quad v\otimes h\mapsto
\sum_{(v)}v_{(0)}\otimes v_{(1)}h.}
Therefore, considering $V$ as an $\sO(L)$-comodule, \eqref{pre06} and \eqref{pre07} imply
a $G$-equivariant isomorphism (where $G$ acts diagonally on the right hand side)
\ga{pre08}{V\Box_{\sO(L)}\sO(G)\cong V\otimes \sO(A).}
In other words, we have $\ind_q\res^q(V)\cong V\otimes\sO(A)$. Note that $\sO(A)$ acts on
$V\otimes\sO(A)$ through the action on the second tensor component.
\subsection{}\label{sect-pre05}
In general, for any $L$-representation $U$ there exists an $\sO(A)$ 
module structure $\mu_U:\sO(A)\otimes \ind_q(U)\longrightarrow \ind_q(U)$
on $\ind_q(U)$, induced from the inclusion of $\sO(A)$ in $\sO(G)$. The action
$\mu_U$ is $G$-equivariant where $G$ acts diagonally on $\sO(A)\otimes\ind_q(U)$.
The functor $\ind_q$ thus factors through a functor to the category $\Mod^{\sO(G)}_{\sO(A)}$
of the so-called ($\sO(G)$-$\sO(A)$)-Hopf modules. 
By definition, an ($\sO(G)$-$\sO(A)$)-Hopf module
is a $k$-vector space $M$ together with a coaction $\rho_M$ of $\sO(G)$ and an action $\mu_M$
of $\sO(A)$ which are compatible in the sense that $\mu_M:M\longrightarrow M\otimes\sO(A)$ is $\sO(G)$-colinear
where $\sO(G)$ coacts diagonally on $M\otimes\sO(A)$. The category
$\Mod^{\sO(G)}_{\sO(A)}$ is a tensor category with respect to the tensor
product over $\sO(A)$. It follows from the various isomorphisms above that
$\ind_q$ defines an equivalence of tensor categories between $\Rep(L)$ and
$\Mod^{\sO(G)}_{\sO(A)}$.
In particular, the isomorphism in \eqref{pre08}
is $\sO(A)$-linear. 
 
 The equivalence mentioned above can be reformulated in the following more geometric
 language: there exists an equivalence between $L$-representations and $G$-equivariant
 vector bundles over $A$. This was pointed out to the author by P. Deligne.
\subsection{}\label{sect-pre06}
A new consequence of the classical facts in \ref{sect-pre01}-\ref{sect-pre06} is the
following, cf. \cite{EsnHai05}.
It follows from the faithful exactness of $\ind_q$ that the canonical homomorphism
$\varepsilon_U:\res^q\ind_q(U)\longrightarrow U$ is surjective. Assume that $U$ has finite dimension over $k$
then we can find a finite dimensional
$G$-subrepresentation $V$ of $\ind_q(U)$ which still maps surjectively on $U$.
Thus $U$ is a quotient of the restriction to $L$ of a finite dimensional representation of
$G$. Consequently, $U$ can also be embedded in to the restriction to $L$ of
a finite dimensional $G$-representation. Thus each finite dimensional representation
$U$ of $L$ can be put (in a non-canonical way) in to a sequence $\res^q(V)
\stackrel\pi\surj U\stackrel\iota\inj \res^q(W)$,
where $V,W$ are finite dimensional $G$-representations. In other words, $U$
is equivalent as an $L$-representation to the image of an $L$-equivariant map
$g:\res^q(V)\longrightarrow\res^q(W)$. Using the equivalence between
$\Rep(L)$ and $\Mod^{\sO(G)}_{\sO(A)}$ and the isomorphism \eqref{pre08} we 
can show that such $g$ are in a 1-1 correspondence with morphisms 
$\bar g:V\otimes\sO(A)\longrightarrow W\otimes \sO(A)$
in $\Mod^{\sO(G)}_{\sO(A)}$. Since $\bar g$ is $\sO(A)$ linear, it is uniquely determined
by a $G$-equivariant map $f:V\longrightarrow W\otimes\sO(A)$.

Let $\sQ$ be the category, whose objects are triples $(V,W,f:V\longrightarrow W\otimes\sO(A))$, where
$V,W$ are  finite dimensional representations of $G$
and $f$ is $G$-equivariant, and whose morphisms are defined in an adequate way. 
Composing $f$ with the morphism $\id\otimes\epsilon:W\otimes\sO(A)\longrightarrow W\otimes k\cong W$
we obtain a map $f_0:V\longrightarrow W$ which is $L$-linear. We define
a functor $\sQ\longrightarrow \Rep(L)$, sending a triple $(V,W,f)$ to the image of $f_0$. It follows
from the discussion of this paragraph that this functor is an equivalence of categories
between $\sQ$ and the category $\Rep_f(L)$ of finite dimensional representations of
$L$.

\section{Groupoids}\label{sect-gr}
\subsection{Groupoids and their representations}\label{sect-gr01}
We refer to \cite[\S1.14]{DeGroth} for the definition of an affine groupoid scheme, which we shall 
call here simply groupoid. Fix a field $k$ and let $K$ be another field containing $k$.
A $k$-groupoid acting upon $\Sp K$ will usually be denoted like $G^K_k$. $G^K_k$ is called
transitive if it acts transitively upon $\Sp K$, which means
that $G^K_k$ is faithfully flat over $\Sp K\times\Sp K$ with respect
to the source and target  map $(s,t):G^K_k\longrightarrow \Sp K\times\Sp K$. 
The category of $K$-representation of
$G^K_k$ is denoted by $\Rep(G^K_k)$, its subcategory of finite dimensional (over $K$) 
representations is denoted by $\Rep_f(G^K_k)$.

For simplicity, we shall omit the subscript $k$ for the product over $\Sp k$. Further, denote $S:=\Sp K$.
\subsection{Homomorphisms of groupoids}\label{sect-gr03}
Assume we are given the following field extensions
$k\subset K_0\subset K$
and transitive groupoids $G^{K}_{k}$ and $A^{K_0}_{k}$.
A $k$-morphism $f:G\longrightarrow A$ is called a groupoid homomorphism if $f$ satisfies
 the following commutative diagram ($S_0:=\Sp K_0$)
\ga{gr06}{\xymatrix{
G^K_k\ar[r]\ar[d]&A^{K_0}_k\ar[d]\\ 
S\times S\ar[r]&S_{0}\times S_0
}}
and is compatible with the groupoid structures on $G$ and $A$
which are related to each other by this diagram.

The \emph{kernel} of $f$ is defined as the fiber product
$L:=G\timesover[A]fe S_0$
where $e:S_0\longrightarrow A$ is  the unite element of $A$.

\subsection{Base change}\label{sect-gr05}
Using the notion of kernel, we define the ``lower" base change for $G_k^K$ 
with respect to $k\subset K_0$ to be the kernel of the homomorphism  
$G^K_k\longrightarrow S\times S\longrightarrow S_0\times S_0$,
and denote it by $G^K_{K_0}$. We notice that
\begin{eqnarray*}
G^K_k\timesover[ S\times S]{}{}( S\times_{S_0} S)
 &=&G^K_k\timesover[ S\times S]{}{}\left(  (S\times S) \timesover[ S_0\times S_0]{}{} S_0\right)=
 G^K_k\timesover[S_0\times S_0]{}{}S_0
 =G^K_{K_0}.
\end{eqnarray*}
Thus $G^K_{K_0}$ can be considered as the restriction of $G^K_k$ to 
$ S_0\times S_0$:
\ga{gr10}{\xymatrix{
\ar@{}[rd]|{\Box} 
G^{K}_{K_0}\ar[r]^{\Delta_{K_0}}\ar[d]&G^K_k\ar[d]\\
 S\times_{S_0}  S\ar[r]& S\times S
}}
In case $K_0=K$, $G^K_K$ is just the usual diagonal subgroup of 
$G^K_k$.

On the other hand, for any extension $K\subset K_1$, Deligne \cite{DeGroth} defines a $k$-groupoid
$G^{K_1}_k$ acting (transitively) upon $ S_1=\Sp(K_1)$:
\ga{gr11}{\xymatrix{
\ar@{}[rd]|{\Box} 
G^{K_1}_{k}\ar[r]\ar[d]&G^K_k\ar[d]\\
 S_1\times S_1\ar[r]& S\times S
}}
which in our context might be called the ``upper" base change.
We notice that the category $\Rep(G^{K_1}_k)$ is equivalent to the category $\Rep(G^K_k)$
(cf. \cite[(3.5.1)]{DeGroth}).

\subsection{The function algebra}\label{sect-gr07}
We refer to \cite[1.14]{DeGroth} or \cite[Appendix]{EsnHai05} for the more detailed desription
 of the function algebra $\sO:=\sO(G^K_k)$ 
of a groupoid $G^K_k$.   Let $m$ denote the multiplication
on $\sO$. Thus $(\sO,m)$ is a faithfully flat $K\otimes K$-algebra by means of the map
$(s,t):s\otimes t:K\otimes K\longrightarrow \sO$ (see \ref{sect-gr01}). The
maps $s,t:K\longrightarrow \sO$ induce two structures of $K$-vector
space on $\sO$, making it a $K$-bimodule. We shall assume that
the left action of $K$ is given by $t$ and the right one by $s$.
Then $\sO$ is a $K$-bimodule coalgebra, where the coproduct
$\Delta:\sO\longrightarrow \sO \otimesover st\sO$ defines the 
product on $G^K_k$, and the counit map $\epsilon:\sO\longrightarrow K$ is the unit element of $G^K_k$.
  Further $\Delta$ and $\epsilon$ are $K\otimes K$-algebra homomorphisms.
 Finally  the antipode map,  $\iota:\sO\longrightarrow \sO$, is induced from the inverse element map 
is an algebra homomorphism and interchanges the actions $s$ and $t$. All these data amount to saying that
$\sO(G)$ is a Hopf algebroid, faithfully flat over $K\otimes K$.

 We adopt Sweedler's notation for the coproduct
on $\sO$:
$\Delta(h)=\sum_{(h)}h_{(1)}{}\otimesover st h_{(2)}.$
We notice that the $K$-linearity of $\Delta$  and $\epsilon$ reads
\ga{gr13}{\Delta(t(\lambda)hs(\mu))=\sum_{(h)}t(\lambda)h_{(1)}{}\otimesover st h_{(2)}s(\mu) ,\quad 
\epsilon(s(\lambda)ht(\mu))=s(\lambda)\epsilon(h)t(\mu),}
where $ h\in \sO$ and $ \lambda,\mu\in K$.
 
The category $\Rep(G^K_k)$ of $G^K_k$ representations over $K$ is the same as
the category of $\sO$-comodules, i.e. of pairs $(V,\rho_V)$ of a $K$-vector space
$V$ and a $K$-linear morphism 
$\rho_V:V\longrightarrow V\otimes_t\sO$
satisfying 
\ga{gr14}{\xymatrix{
V\ar[r]^{\rho_V}\ar[d]_{\rho_V}&V\otimes_t\sO\ar[d]^{\rho_V\otimes\id}\\
V\otimes_t\sO\ar[r]_{\id\otimes\Delta}&V\otimes V\otimes_t\sO
}\qquad
\xymatrix{
V\ar[dr]_{\id}\ar[r]^{\rho_V}&V\otimes_t\sO\ar[d]^{\id\otimes\epsilon}\\
&V
}}
We note that the $K$-linearity of $\rho_V$ means:
$\rho_V(\lambda v)=s(\lambda)\rho_V(v).$
\subsection{The induced representation functor}\label{sect-gr08}
Consider a homomorphism $f:G^{K}_{k}\longrightarrow A^{K_0}_{k}$ of groupoids as in \ref{sect-gr03}.
Let $q:L^{K}_{K_0}\longrightarrow G^{K}_k$ be the kernel of $f$
and $\res^q:\Rep(G^K_k)\longrightarrow \Rep(L^K_{K_0})$ denote the restriction functor.
By definition 
$$\sO(L^K_{K_0})=\sO(G^K_k)\otimes_{\sO(A^{K_0}_k)}K_0,$$
 where $\sO(A^{K_0}_k)$
acts on $\sO(G^K_k)$ through $f^*$, and  $q^*:\sO(G^K_k)\longrightarrow
\sO(L_{K_0}^K)$ is the projection. Define a morphism
\ga{gr16}{\phi:\sO(G^K_k)\otimes_{\sO(A^{K_0}_k)}\sO(G^K_k)
\longrightarrow\sO(L^K_{K_0})\otimesover st\sO(G^K_k)\\
\notag g\otimes h\mapsto \sum_{(g)}q^*(g_{(1)})\otimes g_{(2)}h.} 

The following lemma generalizes \cite[Lem. 6.5]{EsnHai05}.
\begin{lem}\label{lem-gr09}
The morphism $\phi$ in \eqref{gr16} is an isomorphism. It  is $\sO(G^K_k)$-colinear 
 with respect to the right coaction of $\sO(G^K_k)$
on the second tensor component and is $\sO(L^K_{K_0})$-colinear
with respect to the left coaction of $\sO(L^K_{K_0})$ on the first tensor
components.
 \end{lem}
 \proof We give the inverse map.  We first define a map 
$$\bar\psi:\sO(G^K_k){}_s\otimes_t\sO(G^K_k)
 \longrightarrow \sO(G^K_k)\otimes_{\sO(A^{K_0}_k)}\sO(G^K_k) ,
 \quad
\psi(g\otimes h)=
 \sum_{(g)}g_{(1)}\otimes \iota(g_{(2)})h,$$
 where $\iota$ is the antipode of $\sO(G^K_k)$.
 Then we note that this map indeed factors through a map
 $\psi:\sO(L^K_{K_0}){}_s\otimes_t\sO(G^K_k)
 \longrightarrow \sO(G^K_k)\otimes_{\sO(A^{K_0}_k)}\sO(G^K_k)$ which is the inverse to $\phi$.
 The second claim is obvious from the definition of $\phi$.
$\Box$

The functor $\res^q$ possesses a right adjoint which is the induced representation functor,
denoted by $\ind_q:\Rep(L^K_{K_0})\longrightarrow \Rep(G^K_k)$.
For an $\sO(L^K_{K_0})$ comodule $U$ denote  by
$U\Box_{\sO(L^K_{K_0})}\sO(G^K_k)$ the equalizer of the maps
\ga{gr17}{\xymatrix{
\ar[rrd]^{\rho_U\otimes\id}
U\otimes_t\sO(G^K_k)\ar[d]_{\id\otimes\Delta}\\
U\otimes_t\sO(G^K_k) {}_s\otimes_t\sO(G^K_k)
\ar[rr]_{\id\otimes q^*\otimes\id}&&V\otimes\sO(L^K_{K_0}){}
{}_s\otimes_t\sO(G^K_k).
}}

\begin{prop}\label{prop-gr10}  Let $(L=L^K_{K_0},q)$
 be the kernel of $f:G=G^K_k\longrightarrow A=A^{K_0}_k$
as above. Then for $U\in\Rep(L)$, we have
\begin{enumerate}\item[(i)]
 $\ind_q(U)$ is canonically isomorphic to
$U\Box_{\sO(L)}\sO(G)$;
\item[(ii)] the homomorphism
\ga{gr18}{
\left(U\Box_{\sO(L)}\sO(G)\right)\otimes_{\sO(A)}\sO(G)\longrightarrow U\otimes_t\sO(G),\quad
\sum_iu_i\otimes g_i\otimes h\mapsto \sum_iu_i\otimes g_ih
}
is an isomorphism.
\item[(iii)] the algebra
$\sO(G)$ is faithfully flat over its subalgebra $f^*\sO(A)$,
\item[(iv)] the functor $\ind_q:\Rep(L)\longrightarrow \Rep(G)$ is faithfully exact.
\end{enumerate}
\end{prop}
\proof 
(i) We show that the functor $U\mapsto U\Box_{\sO(L)}\sO(G)$ is right adjoint to
the restriction functor $\res^q$, which amounts to
$$\Hom_G(V,U\Box_{\sO(L)}\sO(G)) \cong \Hom_{L}(\res^q(V),U),\quad  V\in\Rep(G), U\in \Rep(L).$$
The map is given by composing a morphism $U\longrightarrow U\Box_{\sO(L)}\sO(G)$ with
the canonical map $\varepsilon_U:U\Box_{\sO(L)}\sO(G)\longrightarrow U$ given by
$v\otimes g\mapsto v\epsilon(g)$, where $\epsilon$ denotes the counit of $\sO(G)$.
The inverse map is given by $f\mapsto (f\otimes\id)\rho_U$. Thus (i) is proved.

(ii) It follows from \eqref{gr16} that the homomorphism
\ga{}{\notag
U\Box_{\sO(L)}\sO(G)\otimes_{\sO(A)}\sO(G)\longrightarrow
U\Box_{\sO(L)}\sO(L)\otimes\sO(G)}
is an isomorphism. Further we have an obvious isomorphism $U\cong U\Box_{\sO(L)}\sO(L)$
given by the coaction $U\longrightarrow U\otimes_t \sO(L)$.

Now (i) and (ii) mean that $\ind_q(-)\otimes_{\sO(A)}\sO(G)\cong (-)\otimes_t\sO(G)$.
Hence  (iii) and (iv) are equivalent, 
for $(-)\otimes_t\sO(G)$ is a faithfully
flat functor.  It thus remains to prove (iii).

(iii)
Let $\bar f:G^K_k\longrightarrow A^K_k=A^{K_0}_k\timesover[S_0\times S_0]{}{}(S\times S)$ be the map obtained from $f$ by means
of \eqref{gr11}. Then the kernel of $\bar f$ is isomorphic to $L^K_K$.
Indeed, we have
$A^K_k\times_{A^{S_0}_k}S_0= S\times_{S_0} S$,
hence
\begin{eqnarray*}
L^K_{S_0}\timesover[ S\times_{S_0} S]{}{}S&=&
( G^K_k\timesover[A^{S_0}_k]{}{}S_0)\timesover[ S\times_{S_0} S]{}{}S=G^K_k\timesover[A^K_k]{}{}
(A^K_k \timesover[A^{S_0}_k]{}{}S_0)\timesover[ S\times_{S_0} S]{}{}S=G^K_k\timesover[A^K_k]{}{}S
 = \Ker\bar f
\end{eqnarray*}
This yields the following commutative diagram with exact lines:
\ga{gr19}{\xymatrix{
&L^{K}_{K}\ar@{^{(}->}[dl]\ar@{^{(}->}[dr]^{\bar q}\ar@{..>}[r]^{q'}&G^K_K\ar@{..>}^\Delta[d] \\ 
L^{K}_{K_0}\ar@{^{(}->}[rr]_q&& G^{K}_{k}\ar[rr]^f\ar[dr]_{\bar f}&& A^{K_0}_{k}\\
&&& A^{K}_{k}\ar@{->>}[ur]
}}

 The morphism $\bar q$ factors as $L^K_K\stackrel{q'}\longrightarrow G^K_K\stackrel{\Delta}\longrightarrow G^K_k$.
 The $K$-group homomorphism $q':L^K_K\longrightarrow G^K_K$ induces a faithful exact functor $\ind_{q'}$, cf. \ref{sect-pre04}.
 Since $G^K_k$ is by assumption faithfully flat over $ S\times  S$,  the isomorphism in 
\eqref{gr18} with $A= S\times  S$ implies that 
the functor $\ind_\Delta$ is faithfully exact. 
 Hence the functor $\ind_{\bar 
q}=\ind_\Delta\circ\ind_{q'}$ is faithfully exact, thus 
 $\sO(G^{K}_{k})$ is faithfully flat over ${\bar f}^*\sO(A^{K}_{k})$ by the equivalence between 
 (iii) and (iv).
 Since $\sO(A^{K}_{k})$ is faithfully flat over $\sO(A^{K_0}_{k})$, we conclude
 that $\sO(G^{K}_{k})$ is faithfully flat over $f^*\sO(A^{K_0}_{k})$.
  $\Box$

  \begin{cor}\label{cor-gr04} $L$ is a $K_0$-groupoid acting transitively upon $ S$.
\end{cor}
\proof This is a direct consequence of Lemma \ref{lem-gr09} and (iii) above. $\Box$

\begin{cor}
(i) Each  finite dimensional representation of $L^K_{K_0}$ is a quotient of the 
restriction to $L$ of a finite dimensional representation of $G^K_k$.

(ii) Each finite dimensional representation of $L^K_{K_0}$ can be embedded into 
the restriction to $L$ of a finite dimensional representation
of $G^K_k$.\end{cor}
\proof
  This follows from Proposition \ref{prop-gr10} (iv) by standard arguments as 
in \cite[Lem. 5.5-5.6]{EsnHai05}. $\Box$
\begin{cor}\label{cor-gr06} 
Let $q:L^K_{K_0}\longrightarrow G^K_k$ be the kernel of a homomorphism
$f:G^K_k\longrightarrow A^{K_0}_k$ of transitive groupoids. For any $V\in \Rep(G^K_k)$, the 
$L^K_{K_0}$-invariant subspace $V^{L^K_{K_0}}$ is a $G^K_k$-subrepresentation of $V$, which
is the restriction along $f$ of a representation of $A^{K_0}_k$.
\end{cor}
\proof  By making the lower base change we see that $L^K_K$ is the kernel of $f_K:G^K_K\longrightarrow A^K_K$,
hence is normal in $G^K_K$ as $K$-groups schemes. Thus $V^{L^K_K}\subset V$
is invariant under the action of $G^K_K$, but this also implies that $V^{L^K_K}$ is invariant under the
action of $G^K_k$. Since $L^K_K$ can also be considered as the kernel of the morphism
$f:G^K_k\longrightarrow A^K_k$, we see that $A^K_k$ acts on $V^{L^K_K}$. 
As we noticed after \eqref{gr11}, $\Rep(A^{K}_k)$ is equivalent to
$\Rep(A^{K_0}_k)$, which means that $V^{L^K_K}$ is indeed a representation of
$G^K_k$ that comes from a representation of $A^{K_0}_k$ by restricting along $f$.
We therefore conclude that $L^K_{K_0}$ acts trivially on this space. Since
$V^{L^K_{K_0}}\subset V^{L^K_K}$, these spaces coincide. $\Box$

\subsection{}\label{sect-gr13} In the situation of Proposition \ref{prop-gr10}, we call a
representation $V$ of $G^K_k$ a $K_0/k$-representation if $V$ is equipped with a $k$-linear
homomorphism $K_0\longrightarrow \text{End}_{G^K_k}(V)$. In other words, there exists a structure of 
$K_0$-vector space on $V$, which is compatible with the $k$-structure and commutes with the
actions of $K$ and $G^K_k$. In the language of comodules,
the above assumption amounts to saying that the comodule map 
$\rho_V:V\longrightarrow V\otimes_t\sO(G^K_k)$ is $K_0$-linear. 
Denote the category of $K_0/k$-representations
by $\Rep(G^K_k)_{K_0/k}$. We notice that the category $\Rep(G^K_k)_{K_0/k}$ 
is a tensor category with respect to the tensor product over $K_0$.

Consider the situation of \ref{sect-gr05}: 
$G^K_{K_0}\xrightarrow{\Delta_{K_0}} G^K_k\longrightarrow S_0\times
S_0$. Then for any representation $W$ over $G^K_{K_0}$,
$\ind_{\Delta_{K_0}}(W)=W\Box_{\sO(G^K_{K_0})}\sO(G^K_k)$ is a $K_0/k$-representation of 
$G^K_k$ with the $K_0$-action induced from the $K_0$-action on $W$.

\begin{lem}\label{lem-gr14}
The functor $\ind_{\Delta_{K_0}}$ induces an equivalence of tensor categories
$\Rep(G^K_{K_0})\longrightarrow \Rep(G^K_k)_{K_0/k}$. In particular, the tensor product over $K_0$
in $\Rep(G^K_k)_{K_0/k}$ is flat.
\end{lem}
\proof To see that $\ind_{\Delta_{K_0}}$ is a tensor functor one has to check that there
is an isomorphism (the cotensor product is over $\sO(G^K_{K_0})$)
\ga{gr20}{(V\Box \sO(G^K_k))\otimesover[K_0]{}{}
(W\Box \sO(G^K_k))\cong 
(V\otimesover[K_0]{}{}W)\Box \sO(G^K_k).
}
We define the map to be 
$\phi:\sum_iv_i\otimes h_i \otimes\sum_jw_j\otimes g_j\mapsto  \sum_{ij}v_i\otimes w_j \otimes  g_ih_j.$
To see that this defines an isomorphism, using the
fact that $\sO(G^K_k)$ is faithfully flat over $\sO(A^{K_0}_k)$, if suffices to show that
\ga{}{\notag (V\Box\sO(G^K_k))\otimesover[K_0]{}{}
(W\Box\sO(G^K_k))\otimesover[\sO(A^{K_0}_k)]{}{}\sO(G^K_k)
\stackrel{\phi\otimes\id}
\longrightarrow(V\otimesover[K_0]{}{}W)\Box\sO(G^K_k)
\otimesover[\sO(A^{K_0}_k)]{}{}\sO(G^K_k)}
is an isomorphism. This last fact follows from the isomorphism in \eqref{gr18}.
$\Box$

\subsection{} \label{sect-gr15}Consider again the situation of \ref{prop-gr10}.
For a representation $W$ of $L^K_{K_0}$, the $K_0$ structure on $W$ yields a $K_0$-structure
on $\ind_{q}(W)=W\Box_{\sO(L^K_{K_0})}\sO(G^K_k)$ making it a $K_0/k$-representation of $G^K_k$.
In particular $f^*\sO(A^{K_0}_k)\cong \ind_q(K)$ is an object of $\Rep(G^K_k)_{K_0/k}$,
moreover, it is an algebra in this tensor category. To simplify the situation, we shall assume that
$f^*$ is injective, i.e., $f$ is a surjective homomorphism of groupoids, and identify $\sO(A^{K_0}_k)$
with its image in $\sO(G^K_k)$.
Denote $\Rep(G^K_k)_{A^{K_0}_k}$ the category of $\sO(A^{K_0}_k)$-modules
in $\Rep(G^K_k)_{K_0/k}$. Then this is a tensor category with respect to the tensor product
over $\sO(A^{K_0}_k)$.
\begin{prop}\label{prop-gr16}With the assumption of \ref{prop-gr10} and 
that $f$ is surjective we have an equivalence of tensor categories $\Rep(L^K_{K_0})\longrightarrow \Rep(G^K_k)_{A^{K_0}_k}$
given by the functor $\ind_q$.
\end{prop}
\proof The proof is similar to that of Lemma \ref{lem-gr14}. $\Box$

\section{Quotient categories}\label{sect-def}
\subsection{Tensor categories and Tannaka duality}\label{sect-def01}
By a {\em rigid tensor category over a field $k$}
 we understand a $k$-linear abelian category $\sC$  equipped with
a symmetric tensor product (see, e.g., \cite[\S1]{DeMil}), such that:
\begin{itemize}\item[(i)]
the endomorphism ring of the unit object (always denoted by $I$) is isomorphic to $k$,
\item[(ii)]  each object is rigid, i.e. possesses a dual object,
and has a decomposition series of finite length.
\end{itemize}
Note that in such a category the hom-sets are
finite dimensional vector spaces
over $k$.

A tensor functor between tensor categories is an additive functor that
preserves the tensor product as well as the symmetry.
A $K$-valued fiber functor of a rigid tensor category $\sC$ over $k$
 ($K\supset k$) is $k$-linear exact tensor functor
from $\sC$ to $\Vect_K$, the category of $K$-vector spaces, its image lies
automatically in the subcategory $\vect_K$ of finite dimensional vector spaces.
If such a fiber functor exists, $\sC$ is called a Tannaka category.
If, more over, $K=k$ then $\sC$ is called neutral Tannaka.

 For example, $\Rep_f(G^K_k)$, where
$G^K_k$ is a transitive groupoid, is a Tannaka category
with the fiber functor being the forgetful functor to $\vect_K$.
 The general Tannaka duality \cite[Thm.~1.12]{DeGroth} establishes
 a 1-1 correspondence between rigid tensor categories over $k$
together with a fiber functor to $\vect_K$ and groupoids acting transitively over $ S$. 

Assume that $\sC$ is rigid over $k$ but not necessarily Tannaka.
Let $\Ind\sC$ be the Ind-category of $\sC$, whose object are
directed systems of objects of $\sC$, and whose hom-sets
are defined as follows:
\ga{d01}{\Hom(X_{i,i\in I}, Y_{j,j\in J}):=\varprojlim_i\varinjlim_j\Hom_{\sT}(X_i,Y_j).}
 The natural inclusion $\sT\hookrightarrow \Ind\sT$
is exact and full. Further, $\sT$ is closed under taking subquotients,
and each object  of $\Ind\sT$ is the limit of its $\sT$-subobjects.  

\begin{defn}\label{def-def02}\rm Let $\sT$ be a rigid tensor 
abelian category over $k$. Let $K\supset k$
be a field extension. 
A normal $K$-quotient of $\sT$ is a pair $(\sQ,\fq:\sT\longrightarrow\sQ)$
consisting of a $K$-linear rigid tensor category $\sQ$ and 
 $k$-linear exact tensor functor $\fq$ (the $k$-linear structure over
$\sQ$ is induced from the inclusion $k\subset K$), such that:
\begin{itemize}\item[(i)]for an object $X$ of $\sT$,
the largest trivial subobject of $\fq(X)$  is isomorphic to the image
under $\fq$ of a subobject of $X$;
\item[(ii)] each object of $\sQ$ is isomorphic to a subobject of the image 
of an object from $\sT$ as well as a quotient of the image of an object from $\sT$.
\end{itemize}
\end{defn}
For convenience  we 
shall omit the term ``normal'' in the rest of the work.
Our notion of quotient category in case 
$K=k$ and $\sT$ is a Tannaka category is equivalent to Milne's notion of normal quotient 
 \cite{milne2005}. 
Indeed, assume that there exists a fiber functor $\bar\omega:\sT\longrightarrow \vect_{K}$
extending the fiber functor $\omega$ and
$\sT\cong\Rep_f(G^K_k)$, $\sS\cong\Rep_f(A^k_k)$ for some
groupoid scheme $G^K_k$ and group scheme $A_k^k$.
According to Corollary \ref{cor-gr06}, Proposition \ref{prop-gr10},
$\sQ:=\Rep_f(L^K_k)$ is a $k$-quotient of $\sT$ by $\sS$.

Given a $K$-quotient $(\sQ,\fq)$ of $\sT$, let $\sS$ denote the full
subcategory of $\sT$ consisting of those objects of $\sT$
whose images in $\sQ$ are trivial (i.e. isomorphic to a direct
sum of the unit object). It is easy to see that
$\sS$ is a tensor subcategory and is closed under taking subquotients.
We shall call $\sS$ the \emph{invariant subcategory} with respect to
the quotient $(\sQ,\fq)$.

\begin{lem}\label{lem-def03} \cite[\S 2]{milne2005}
Let $(\sQ,\fq)$ be a $K$-quotient of $\sT$ and $\sS$ be the
invariant subcategory of $\sT$. Then $\sS$ is a Tannaka
category with a fiber functor to $\vect_K$.
\end{lem}
{\it Proof.}   
The full subcategory of $\sQ$ of trivial subobjects is 
equivalent to $\vect_K$.
The fiber functor is given by
$\sS\ni X\mapsto \fq(X)\mapsto \Hom_{\sQ}(I,\fq(X))$
where $I$ denotes the unit object in $\sQ$. Since $\fq(X)$ is trivial in
$\sQ$, this functor is a fiber functor. Hence $\sS$ is a Tannaka category. \hfill
$\Box$
\section{Quotient category by a neutral Tannaka subcategory}\label{sect-q}

Let $\sT$ be a rigid tensor abelian category over a field $k$.
Let $\fq:\sT\longrightarrow\sQ$ be a $k$-quotient of $\sT$ and
$\sS$ be the invariant category as in Definition \ref{def-def02} (with $K=k$). 
Then  
$\sS$ is a neutral Tannaka category with fiber functor given
by $\omega(S)\cong \Hom_{\sQ}( I,\fq S)$, cf. \ref{lem-def03}. Let us consider
the category $\vect_k$ as a full subcategory of
$\sQ$ by identifying a vector space $V$ with $V\otimes I$
in $\sQ$. Then we can consider the above 
fiber functor as the restriction of $\fq$ to $\sS$. In other words,
 we have the following functorial isomorphism
\ga{q01}{\omega(S)\otimes_k I= \Hom_\sQ(I,\fq(S))\otimes_kI\cong \fq(S),\quad S\in \sS,}
by means of which we shall identify $\omega(S)$ with $\fq(S)$ for 
$S\in\sS$.
Our aim is to characterize $\sQ$ in terms of $\sS$ and deduce (with some technical
condition) the uniqueness of $\sQ$ with given $\sS$. We also show the existence
of $\sQ$ in case $\sS$ is a finite tensor category, i.e. its Tannaka group is finite.

\subsection{The algebra $\sO$}\label{sect-q0a}
By means of the fiber functor $\omega$,  $\sS$ is equivalent to the category
$\Rep_f(G(\sS))$ of finite dimensional $k$-representation of $G(\sS)$, where
$G(\sS)$ is an affine $k$-group scheme. Further $\omega$ extends to an equivalence between
$\Ind\sS$ and $\Rep(G(\sS))$. 
Let $\sO $ denote the function algebra of  the 
group $G(\sS)$. It is a $k$-Hopf algebra. 
By means of the right regular action of $G(\sS)$ on $\sO $,
that is, consider $\sO $ as a right comodule on itself by means
of the coproduct map,  $\sO $  can be considered 
as an object in $\Ind\sS$. Since the
unit map $u:k\longrightarrow \sO$  and the multiplication map
$m:\sO\otimes_k\sO\longrightarrow \sO$ of $\sO $ are compatible with
the coproduct and the counit maps, they are also morphisms in $\Ind\sS$, hence
$(\sO,m,u)$ is an algebra in $\Ind\sS$ (it is not a Hopf algebra since the
coproduct and the counit are not morphisms in $\Ind\sS$).
\subsection{Convention}\label{sect-q0b}
For convenience we shall use the same notation for denoting $\sO $
as a $k$-vector space or as an object in $\sS$ as well as an object of
$\sQ$ when we consider $\omega$ as the restriction of $\fq$ to $\sS$.

 \subsection{The largest $\sS$-subobject}
 For an object $X$ of $\sT$, let $X_\sS$ denote the largest $\sS$-subobject of $X$
(which exists by the assumption \ref{sect-def01} (ii)). 
 Since  $\sS$ is closed under taking subquotients, we have the equality
  \ga{q02}{ \Hom_\sT(S,X)= \Hom_\sS(S,X_\sS),\quad S\in\sS, X\in\sT.}
Thus we have the functor
 $(-)_\sS:\sT\longrightarrow \sS$ whose definition on hom-sets is just the restriction of
 morphisms
$\Hom_{\sT}(X,Y)\mapsto\Hom_{\sT}(X_{\sS}, Y)=\Hom_{\sT}(X_{\sS}, Y_{\sS})$.
  Equation (\ref{q02}) also shows that this functor is right adjoint to the 
  inclusion functor $\sS\hookrightarrow\sT$.
 The functor $(-)_\sS$ is canonically extended to a functor
 $\Ind\sT\longrightarrow \Ind\sS$, (denoted by the same symbol), 
which is also the right adjoint
 to the inclusion functor $\Ind\sS\longrightarrow\Ind\sT$.
 
 Let $X^\vee$ denote the dual object to $X$. The for any $S\in \sS$, we have
 $$\Hom_\sT(X,S)\cong\Hom_{\sT}(S^\vee,X^\vee)\cong
 \Hom_{\sS}(S^\vee,(X^\vee)_\sS)\cong \Hom_{\sS}((X^\vee)_\sS{}^\vee,S).$$
 Since $X_\sS\longrightarrow X$ is mono, $X^\vee\longrightarrow (X^\vee)_\sS{}^\vee$
 is epi. Thus the largest $\sS$-quotient of $X$ is isomorphic to
 $(X^\vee)_\sS{}^\vee$.
 \begin{lem} \label{lem-q01} There is a functorial isomorphism
\ga{q03}{\Hom_{\Ind\sT}(I,X\otimes\sO )\stackrel\cong\longrightarrow
\omega(X_{\sS}),\quad  X\in \sT.} 
\end{lem}
\proof  If $X=S$ is an object of $S$ we have the following  isomorphism (see \eqref{pre07})
\ga{q04}{S\otimes\sO \cong \omega(S)\otimes_k \sO .} 

For $X$ arbitrary, $\Ind\sT$-morphisms $I\longrightarrow X\otimes\sO $ are in 1-1
correspondence with $\Ind\sT$-morphisms $X^\vee\longrightarrow \sO $, which, are in
1-1 correspondence with morphisms from the largest $\sS$-quotient
of $X^\vee$ to $\sO $, since $\sO $ is an object of $\Ind\sS$.
Notice that the largest $\sS$-quotient of $X^\vee$ is isomorphic to
$(X_\sS)^\vee$, as for  $X\in \sT$  one has $(X^\vee)^\vee \cong X$. Thus we have
$$
\Hom_{\Ind\sT}(I,X\otimes \sO )\cong\Hom_{\Ind\sT}(X^\vee,\sO )
\cong\Hom_{\Ind\sT}((X_\sS)^\vee,\sO )\cong\omega(X_\sS).
$$ 
Therefore (\ref{q03}) is proved. The functoriality of (\ref{q03}) is obvious
since any morphism $f:X\longrightarrow Y$ restricts to a morphism $f:X_\sS\longrightarrow Y_\sS$.
 $\Box$
\bigskip
\begin{prop}\label{prop-q02} (i) There exists a functorial isomorphism
\ga{q06}{\Hom_\sQ(\fq(X),\fq(Y))\cong\Hom_{\Ind\sT}(X,Y\otimes \sO )\\
\notag f\mapsto \bar f}
satisfying the following diagram
\ga{q05}{\xymatrix{ &\fq(Y)\otimes \sO\ar[d]^{\id\otimes\epsilon}\\
\fq(X)\ar[ru]^{\fq(\bar f)}\ar[r]_f&\fq(Y)}}
 
(ii)  For morphisms $f:\fq(X)\longrightarrow \fq(Y)$,
 $g:\fq(Y)\longrightarrow \fq(Z)$, the morphism $\overline{gf}:
X\longrightarrow Z\otimes \sO$ is given by
$$ X\stackrel {\bar f} \longrightarrow Y\otimes \sO \stackrel{\bar g\otimes \id}\longrightarrow Z\otimes \sO \otimes
 \sO \stackrel{\id\otimes m}\longrightarrow Z\otimes\sO. $$
 
 (iii) For $f_i:\fq(X_i)\longrightarrow \fq(Y_i)$, $i=0,1$  the morphism $\overline{f_0\otimes f_1}$ is given by
 $$\xymatrix{X_0\otimes X_1\ar^{\overline{f_0\otimes f_1}}[r]
 \ar_{\bar f_0\otimes \bar f_1}[d]&
  Y_0\otimes Y_1\otimes\sO \\ 
Y_0\otimes\sO \otimes Y_1\otimes \sO \ar_{\tau_{(23)}}[r]&
Y_0\otimes Y_1\otimes \sO \otimes\sO  
\ar^m[u]
} $$
where the map $\tau_{(23)}$ interchanges the second and the third 
tensor components.
\end{prop}
\proof  (i). 
The isomorphism in \eqref{q06} is defined as follows:
\ga{}{\notag\xymatrix{
\Hom_{\sQ}(\fq(X), \fq(Y)) \ar[d]_\cong \ar@{..>}[rr] && \Hom_{\Ind\sT}(X,Y\otimes \sO )\\
\Hom_{\sQ}(I,\fq(X^\vee\otimes Y)) \ar[rr]_{\eqref{q02},\eqref{q03}\quad}&& \Hom_{\Ind\sT}(I,X^\vee\otimes Y\otimes\sO )
\ar[u]_\cong
}}
It follows from the definition that \eqref{q05} needs only to be checked for
$X=I$ and hence one can replace $Y$ by $Y_{\sS}$. But in this case the claim
follows from \eqref{q04} and \eqref{pre07}.

(ii).  Since $\fq$ is a faithful functor, it suffices to check that the outer diagram below commutes.
 $$\xymatrix{
 \fq(X) \ar^{\fq(\bar f)\quad}[r]\ar_{f}[dr]
&\fq(Y)\otimes\sO  
\ar^{\id\otimes\epsilon}[d]
\ar^{\fq(\bar g)\otimes\id\hspace{3ex} }[r]
&\fq(Z)\otimes\sO \otimes\sO \ar^{\quad m}[r]
\ar_{\id\otimes\epsilon\otimes\epsilon}[d]
&\fq(Z)\otimes\sO \ar^{\id\otimes\epsilon}[dl]\\
&\fq(Y)\ar_{g}[r]& \fq(Z)
}$$
The commutativity of the first triangle and the middle square follow
from (\ref{q05}) and the commutativity of the right triangle is by the multiplicativity
of the counit map $\epsilon$.

The proof for (iii) is similar and will be omitted.
 $\Box$
 
\subsection{The adjoint functor to $\fq$}\label{sect-q03}
 The quotient functor $\fq$ extends to a functor from 
 $\Ind\sT\longrightarrow \Ind\sQ$, and will be denoted by the same symbol.
 This functor possesses a right adjoint, denoted by $\fp$
 \ga{q08}{ \Hom_{\Ind\sQ}(\fq(  X),  U)\cong 
 \Hom_{\Ind\sT}(  X,\fp(  U)),\quad \forall  X
 \in\Ind\sT,   U\in\Ind\sQ.}
 Indeed, for $U\in \sQ$, the above isomorphism determines $\fp(U)$.
 Further the definition \eqref{d01} of the hom-sets of $\Ind\sT$ shows that
  $\fp(\varinjlim_iU_i)=\varinjlim_i\fp(U_i).$ 

Our next aim is to show that
  there is a natural $\sO$-module structure on $\fp(U)$ for $U$ from $\sQ$.
We first recall that there is a natural morphism $\varepsilon$ associated to
the adjoint pair $(\fq,\fp)$.
  For $X=\fp( {U})$ in \eqref{q08}, the identity of $\fp(  U)$ corresponds to a map 
 $\varepsilon_{ {U}}:\fq\fp( {U})\longrightarrow  {U}$, and \eqref{q08} can be given 
 by composing a morphism on the right  hand side with  $\varepsilon_{ {U}}$:
\ga{q07}{\xymatrix{
  &\fq\fp(U)\ar^{\varepsilon_U}[d]\\
  \ar^{\fq(\bar f)}[ur]  \fq(X)\ar_f[r]&U}  
}
  \begin{lem}\label{lem-q04} For an object $X$ of $\Ind\sT$, the object
 $\fp\fq(X)$ in $\Ind\sT$  
 is canonically isomorphic to $X\otimes\sO $ and the map 
 $$\fq( X)\otimes\sO \cong\fq\fp\fq(X)\stackrel{\varepsilon_{\fq(X)}}\longrightarrow \fq(X)$$
 is given by $\id_X\otimes\epsilon$,
where $\epsilon$ is the counit of $\sO $. With respect to the isomorphism
$\fp\fq(X)\cong X\otimes \sO$, (\ref{q08}) for $U=\fq(Y)$ reduces to (\ref{q06}).
Consequently, for a morphism $f:\fq(X)\longrightarrow \fq(Y)$, $\fp(f)$ is given by
 \ga{q11}{X\otimes\sO \stackrel{\bar f\otimes\id}\longrightarrow Y\otimes\sO \otimes\sO 
 \stackrel{\id\otimes m}\longrightarrow Y\otimes\sO.}
\end{lem}

\proof  First assume $X$ is in $\sT$.
 According to  \eqref{q06} and (\ref{q08}) we have, for $X,Y\in\sT$,
 \ga{q09}{\Hom_{\Ind\sT}(X,\fp(Y)) \cong
 \Hom_{\Ind\sT}(X,Y\otimes\sO )}
which implies that $\fp(Y)$ is canonically isomorphism to $Y\otimes \sO $.
The second claim follows from \eqref{q05} and \eqref{q07}.

 For the general case ($X$ is in $\Ind\sT$) we notice that the tensor product in $\Ind\sT$ commutes
 with direct limits, hence
$$\fp\fq(\varinjlim_iX_i)\cong \varinjlim_i (X_i\otimes \sO)\cong 
(\varinjlim_i X_i)\otimes \sO.$$

To prove the last claim we notice that $\overline{\id_{\fq(X)}}=\id_X\otimes U:
X\longrightarrow X\otimes\sO$. Indeed, this map satisfies \eqref{q07} with $U=\fq(X)$,
since $\varepsilon_{\fq(X)}=\id_X\otimes\epsilon$. On the other hand, it follows
from \eqref{q08} that for $f:\fq(X)\longrightarrow\fq(Y)$, $\fp(f)$ is determined
by the diagram
 $$\xymatrix{X\ar[r]^{\overline{\id_{\fq(X)}}} \ar[dr]_{\bar{f}}&\fp \fq(X)\ar[d]^{\fp(f)}\\
& \fp\fq(Y) 
 } $$
The morphism in (\ref{q11})
 satisfies this property, hence is equal to $\fp(f)$.  
  $\Box$
  
Since $\sO$ is a commutative algebra in $\Ind\sT$, we can  consider the
category $\Mod_\sO$ of $\sO$-modules of $\Ind\sT$, which is an abelian category equipped with a tensor product over $\sO$.
\begin{defn} \label{def-q08}\rm 
Let $(\sS,\omega:\sS\longrightarrow\vect_k)$ be a neutral Tannaka subcategory of a rigid tensor abelian category $\sT$ which is closed under
taking subquotient objects. Denote by $\sO$ the function algebra over its Tannaka group,
viewed as an object in $\Ind\sS\subset\Ind\sT$. The category $\sT$ is said to 
be \emph{flat} over $\sS$ if the tensor product in $\Mod_\sO$ is exact.
\end{defn}
 \begin{thm} \label{thm-q09}Let $(\sQ,\fq)$ be a $k$-quotient category
 of a rigid tensor abelian category $\sT$ over $k$ and denote by $\sS$ the corresponding
 invariant category. Let $\sO$ denote the function algebra over
 the Tannaka group of $\sS$. Let $\fp$ be the right adjoint to the
 quotient functor $\fq$. Then:
 \begin{enumerate}\item For any object $U\in\sQ$, there exists a morphism
 $\mu_U:\fp(U)\otimes\sO\longrightarrow \fp(U)$ making $\fp(U)$ an $\sO$-module.
 \item Assume that $\sT$ is flat over $\sS$. Then
$\fp$ is a tensor functor from $\sQ$ to the category
 of $\sO$-modules with the tensor product being the tensor product over $\sO$.
 Consequently, $\fp$ is exact. \end{enumerate}
 \end{thm}
 \proof 
Recall that we have $\fp\fq(X)=X\otimes\sO$ and that we identify, by means of $\omega$,
$\sO$ and its image in $\sQ$.
Using (\ref{q08}), we define $\mu_U$ as the unique morphims in $\Ind\sT$ making
 the following diagram commutative (cf. \eqref{q07})
 $$\xymatrix{&\fq\fp(U)\ar[d]_{\varepsilon_U}\\
 \fq\fp(U)\otimes \sO\ar[ur]^{\fq(\mu_U)}\ar[r]_{\qquad\qquad 
\varepsilon_U\otimes\epsilon}&U
 }$$
 In other words, $\mu_U:=\overline{\varepsilon_U\otimes\epsilon}$.
 This definition is functorial hence the action of $\mu$ commutes with morphisms
 in $\sQ$. The associativity of this action can also be checked by this method. Thus
 $\fp$ factors through a functor to $\Mod$-$\sO$, denoted by the same notation.

For  $U=\fq(X)$,
 $\fp \fq(X)=X\otimes\sO$, we see that the multiplication map on $\sO$ makes
 the following diagram commutative
 $$\xymatrix{&\fq\fp(\fq(X))= \fq(X)\otimes\sO\ar[d]_{\varepsilon_{\fq(X)}}\\
 \fq(X) \otimes\sO\otimes\sO\ar[ur]^{ \id\otimes m}\ar[r]_{\qquad\qquad \varepsilon_{\fq(X)}\otimes\epsilon}&\fq(X)
 }$$
 Thus, the action of $\sO$ on $\fp\fq(X)= X\otimes\sO$ is induced from the action
 of $\sO$ on itself. 

 Assume now that $\sT$ is flat over $\sS$.
  We want to show that $\fp$ is a monoidal functor from $\sQ$ to $\sO$-modules,
  that is to show the existence of a functorial isomorphism
  $$\fp(U_1)\otimes_{\sO}\fp(U_2)\cong\fp(U_1\otimes U_2).$$
 This is so for objects of the form $\fq(X)$, according to Proposition \ref{prop-q02} (iii).
 For an arbitrary object $U$ of $\sQ$, there exist objects $X,Y$ in $\sT$ and a morphism
 $f:\fq(X)\longrightarrow \fq(Y)$ such that $U=\ker f$.  
Since the functor $\fp$ is left exact (being a right adjoint functor), $\fp(U)$ is isomorphic
to the kernel of $\fp(f)$. Now, for $f_i:\fq(X_i)\longrightarrow \fq(Y_i)$ we have $\fp(f_1\otimes f_2)
=\fp(f_1)\otimes_{\sO}\fp(f_2)$. Since for $U_i:=\ker f_i$ are flat over $\sO$ by assumption, 
we have an $\sO$-linear isomorphism
$$\fp(U_1)\otimes_{\sO}\fp(U_2)\cong\fp(U_1\otimes U_2).$$
Thus $\fp$ is a tensor functor to $\sO$-modules.
Since $\sQ$ is rigid with the endomorphism ring of the unit object isomorphic
to $k$ and $\Mod_\sO$ has an exact tensor product,
we conclude that $\fp$ is an exact functor, cf. \cite[2.10]{DeGroth}. $\Box$

\subsection{A description of the quotient}\label{sect-q10}
 We define a category $\sP$, whose objects are triples
$$(X,Y,f:X\longrightarrow Y\otimes\sO),$$
where $X,Y$ are objects of $\sT$ and $f$ is a morphism in $\Ind\sT$.
 The morphisms $f:X\longrightarrow Y\otimes \sO$ defines an $\sO$-linear morphism
 $\hat f$
\ga{q12}{\hat f :X \otimes\sO \stackrel{f \otimes \id}\longrightarrow
Y \otimes \sO \otimes \sO \stackrel {\id\otimes m}\longrightarrow Y \otimes \sO .}
The image of $\hat f$ is thus an $\sO$-module.

We  define morphisms in $\sP$ between
 $(X_i,Y_i,f_i)$, $i=0,1,2$ as $\sO$-linear morphisms $\phi:\im \hat f_0\longrightarrow \im \hat f_1$.
Thus we obtain a category $\sP$ which is $k$-linear and additive. Further, the
direct sum of objects in $\sP$ exists:
\ga{q13}{(X_0,Y_0,f_0)\oplus (X_1,Y_1,f_1):=
(X_0\oplus X_1,Y_0\oplus Y_1,f_0\oplus f_1).}

The tensor structure on $\sP$ is defined as follows.
\ga{q14}{(X_0,Y_0,f_0)\otimes (X_1,Y_1,f_1):=
(X_0\otimes X_1,Y_0\otimes Y_1,(\id\otimes m)(f_0\otimes f_1)).}
The unit object is $I=(I,I,u:I\longrightarrow \sO)$.
Finally, we define a functor $\fq':\sT\longrightarrow \sP$ sending
$X$ to $\fq'(X)=(X,X,\id\otimes u:X\longrightarrow X\otimes \sO)$ and sending $f:X\longrightarrow Y$ to
$f\otimes\id_\sO$.
 It is easy to see that $\fq'$ is a $k$-linear,
additive tensor functor.

\begin{prop}\label{prop-q11}
Let $\sT\stackrel\fq\longrightarrow \sQ$ be $k$-quotient and $\sS$ be the
corresponding invariant subcategory. Assume that $\sT$ is flat over $\sS$.
Then the category $\sP$ constructed above is equivalent to $\sQ$.
\end{prop}
\proof 
We construct a functor $\ff:\sP\longrightarrow \sQ$ such that $\ff\fq'=\fq$.
 For any $X,Y\in \sT$,  according to \eqref{q06}, morphisms $f:X\longrightarrow Y\otimes\sO $ in $\sT$
are in 1-1 correspondence with morphisms $f_{\fq}:\fq(X)\longrightarrow \fq(Y)$ in $\sQ$.
This allows us to define the image under $\ff$ of an object
$U=(X,Y,f:X\longrightarrow Y\otimes\sO )\in \sQ'$ as the image of $f_{\fq}$ in $\sQ'$.

According to Theorem \ref{thm-q09}, $\fp$ is exact. Hence for any morphism
$f:\fq(X)\longrightarrow \fq(Y)$ in $\sQ$ with $U=\im f$ one has
$\fp(U)\cong \im\fp(f).$
Consequently, there is a 1-1 correspondence between morphisms
$U\longrightarrow V$ in $\sQ$ and morphisms between
$\fp(U)\longrightarrow \fp(V)$ in $\sT$ which are $\sO $-linear. This allows
us to define $F$ on morphisms and to check that $F$ is a fully faithful
functor.
  Further it is easy to see that the image of $F$ is essential in
$\sQ$.  Thus $F$ is an equivalence. 
$\Box$

As a consequence of the above theorem and Proposition \ref{prop-gr16} we have
the following.
\begin{cor}\label{cor-q12} Let $f:G^K_k\longrightarrow A^k_k$ be a homomorphism
of transitive groupoids and $L^K_k$ be the kernel of $f$. Then $\Rep_f(L^K_k)$
is equivalent to the category whose objects are triples $(U,V,f:U\longrightarrow V\otimes\sO(A))$,
where $U,V\in\Rep_f(G)$, $f$ is $G$-equivariant ($G$ acts diagonally on $V\otimes\sO(A))$.
\end{cor}
\subsection{An alternative description} \label{sect-q15}
There is an alternative description of objects of $\sQ$, proposed
by P.~Deligne. Notice that morphisms $f: X\longrightarrow Y\otimes\sO$ are in 1-1 correspondence
with morphisms $f^\sharp:I\longrightarrow X^\vee\otimes Y\otimes\sO$. 
As noticed in \ref{lem-q01}, such a
morphism $f^\sharp$ corresponds to an element of $\omega((X^\vee\otimes Y)_\sS)$.
Thus objects of $\sQ$ can be characterized as triples
$(X,Y,f\in\omega((X^\vee\otimes Y )_\sS))$.

\subsection{The case $\sS$ is an \'etale finite tensor category}\label{sect-q16}
The (neutral) Tannaka category $\sS$ is said to be finite if its Tannaka group is finite, this means
the algebra $\sO$ lies in $\sS$ and hence in $\sT$.
By abusing language, we say that $\sS$ is \'etale finite if its Tannaka group is finite and reduced,
this amount to saying that $\sO$, as a $k$-algebra, is a separable algebra, i.e. $\sO$ is
projective as an $\sO$-bimodule.

 Our next aim is to show that if $\sS$ is  \'etale finite  then $\sT$ is flat
over $\sS$ and, moreover, the quotient of $\sT$ by $\sS$ exists.
Consider $\sO$ as a $k$-algebra. 
If $\sO$ is separable then the multiplication map $\sO\otimes \sO\xrightarrow m \sO$,
considered as an $\sO$-bimodule map, has a section. Consequently, for any $\sO$-module $M$, taking the
tensor product over $\sO$ of $M$ with the split exact sequence $0\to J\to \sO\otimes \sO\xrightarrow{m} \sO\to 0$
we obtain a split exact sequence of $\sO$-modules 
$$0\to J\otimes_{\sO}M\to \sO\otimes_kM\xrightarrow{\mu_M} M
\to 0$$
where $\mu_M$ denote the multiplication of $\sO$ on $M$. This implies that $M$, being a direct summand
of the flat $\sO$-module $\sO\otimes_kM$, is flat as well. This construction
will be use to prove the following theorem.
\begin{thm}\label{thm-q18} Let $(\sS,\omega)$ be a neutral Tannaka subcategory of 
$\sT$ which is closed under taking subquotients. Assume that $\sS$ is \'etale finite.
Then the category of $\sO$-modules
in $\sT$ is rigid. Consequently $\sT$ is flat over $\sS$ and the quotient of
$\sT$ by $\sS$ exists.
\end{thm}
\proof 
To apply the above method to our case, where $\sO$ is also considered as a comodule on itselft,
we have to find a splitting of $m:\sO\otimes\sO\to \sO$, which is $\sO$-colinear.
This is done as follows. Let $x$ be a non-zero integral element in $\sO$, that is, $x$ satisfies the 
condition $ax=\varepsilon(a)x$ for any $a$ in $\sO$. It is due to
a theorem of Sweedler \cite[Theorem~5.1.2]{sweedler} for finite dimensional (not necessarily commutative)
Hopf algebras that $x$ exists uniquely up to a scalar constant. Let $e_i$ be a $k$-basis of
$\sO$ and $e^i$ be elements of $\sO$ such that
$$\Delta(x)=\sum_i e_i\otimes e^i.$$
We claim that there is an $\sO$-bimodule map
 $s:\sO\to \sO\otimes \sO$, with $s(1)=\sum_i e_i\otimes \iota(e^i)$
($\iota$ denotes the antipode of $\sO$), which
splits the multiplication map $m$ up to some non-zero scalar constant and is $\sO$-colinear.

In fact, one can assume without lost of generality that $k$ is algebraically closed. In this case,
by the \'etale assumption on $\sO$, there is a finite group $G$ such that $\sO=k[G]^*$,
i.e., it has a basis $x^g, g\in G$ with multiplication, comultiplication, unit, counit given as follows:
$$x^gx^h=\delta^g_hx^h, \textstyle\quad u(1)=\sum_g x^g$$
$$\Delta(x^g)=\sum_h x^h\otimes x^{h^{-1}g},\quad \epsilon(x^g)=\delta^g_e,$$
($e$ denotes the unit element of $G$ and $\delta$ denotes the Kronecker symbol) and the antipode is 
$$\iota(e^g)=e^{g^{-1}}.$$
Now the integral element of $\sO$ is (up to a constant) $x=x^e$. Thus the map $s$ is
$$s(1)=\sum_ge^g\otimes e^g.$$
It is now strightfoward to check that $s$ defines an $\sO$-bimodule map
 (i.e. $s(1)$ is central in $\sO\otimes \sO$),
which splits $m$ and is $\sO$-colinear.

Let now $(M,\mu_M)$ be an $\sO$-module in $\sT$. 
Consider $\sO\otimes M$ as an $\sO$-module by the multiplication on $\sO$.
Then the map $\mu_M:\sO\otimes M\to M$ is $\sO$-linear and moreover
can be obtained by tensoring the exact sequence $0\to J\to\sO\otimes\sO\xrightarrow m\sO$
with $M$ over $\sO$. Since the last exact sequence split as a sequence of $\sO$-bimodules,
$\mu_M$ has a section which is $\sO$-linear. Thus $M$ is a direct summand
 of $\sO\otimes M$ as an $\sO$-submodule. Since $\sO\otimes M$ 
(with the $\sO$-action given by the multipliacation on $\sO$) is flat and rigid as
 an $\sO$-module, so is $M$.  
 
 Thus we have shown that the category of $\sO$-modules in $\sT$ is a $k$-linear, rigid
 tensor category. It is now clear that the functor $\fq: X\mapsto \sO\otimes X$ makes
 this category a quotient tensor category of $\sT$ by the Tannaka subcategory $(\sS,\omega)$.
$\Box$

\section{Base change for tensor category}\label{sect-b}
\subsection{}\label{sect-b01}
Let $\sT$ be a tensor category over $k$. Assume that for a field extension
$k\subset K$, a $K$-quotient $\sT_{(K)}$ of $\sT$ exists such that the invariant
subcategory is the trivial subcategory of $\sT$, i.e. equivalent to $\vect_k$.
$\sT_{(K)}$ is called the tensor category over $K$ obtained from $\sT$ by
base change. Denote the quotient functor by $-_{(K)}:X\mapsto X_{(K)}$.
Thus, the largest trivial subobject $X_{(K)}{}^{\rm triv}$ of $X_{(K)}$ in $\sT_{(K)}$
is isomorphic to the image of the largest trivial object $X^{\rm triv}$ of $X$ in
$\sT$. Hence
\begin{eqnarray*}
\Hom_{\sT_{(K)}}(I,X_{(K)})= \Hom_{\sT_{(K)}}(I,X_{(K)}{}^{\rm triv}) 
=\Hom_{\sT}(I,X^{\rm triv})\otimes_k K 
= \Hom_{\sT}(I,X)\otimes_kK.
\end{eqnarray*}
Consequently
\ga{b01}{\Hom_{\sT_{(K)}}(X_{(K)},Y_{(K)})=\Hom_{\Ind\sT}(X,Y)\otimes_kK.}

On the other hand, we have an isomorphism
\ga{b02}{\Hom_{\sT}(X,Y)\otimes_kK\cong 
\Hom_{\Ind\sT}(X,Y\otimes_kK),}
 given explicitly as follows.
Fix a basis $\{e_i,i\in I\}$ of $K$ over $k$. An element $f\in \Hom_{\sT}(X,Y)\otimes_kK$,
represented as $f=\sum_{i\in I_0}f_i\otimes e_i$, for a certain finite subset $I_0\subset I$,
is mapped to the morphism
$$\bar f:X\stackrel\Delta\longrightarrow \bigoplus_{i\in I_0}X_i\stackrel{\oplus_{i}f_i}\longrightarrow
\bigoplus_{i\in I_0}Y_i\inj Y\otimes_kK,$$
where $X_i$ (resp. $Y_i$) are copies of $X$ (resp. $Y$) and
 the last inclusion in given by the chosen basis of $K$.
 The bijectivity of the correspondence follows from the definition
 of morphism in $\Ind\sT$. It follows from \eqref{b01} and \eqref{b02}
 that
 \ga{b03}{\Hom_{\sT_{(K)}}(X_{(K)},Y_{(K)})=\Hom_{\Ind\sT}(X,Y\otimes_kK),\quad f\mapsto
 \bar f.}

Extend $-_{(K)}$ to a functor $\Ind\sT\longrightarrow\Ind\sT_{(K)}$ and 
let ${-}^K:\Ind\sT_{(K)}\longrightarrow \Ind\sT$ denote the right adjoint functor to $-_{(K)}$. 
Claims similar to those of Proposition \ref{prop-q02}
and of Lemma \ref{lem-q04} hold true.  We collect them in
the following lemma.
\begin{lem}\label{lem-b02}
(i) Through the isomorphism \eqref{b03}, the composition
$X_{(K)}\stackrel{f}\longrightarrow Y_{(K)}\stackrel{g}\longrightarrow Z_{(K)}$ corresponds to
$X\stackrel{\bar f}\longrightarrow Y\otimes_k K\stackrel{\bar g\otimes \id}\longrightarrow Z\otimes_kK\otimes_kK
\stackrel{\id\otimes m}\longrightarrow Z\otimes_kK$, (cf. \ref{prop-q02} (ii)), 
and the tensor product $f_0\otimes f_1$ of $f_i:X_{i(K)}\longrightarrow Y_{i(K)}, i=0,1$ corresponds to
$ (\id\otimes m)\tau_{(23)}(\bar f_0\otimes \bar f_1)$ (cf. \ref{prop-q02} (iii)).

(ii) For an object $X$ of $\sT$ holds: $(X_{(K)})^K\cong X\otimes_kK$,
and the map $\varepsilon_X$ is given
by $\varepsilon_X:X_{(K)}\otimes_kK\longrightarrow X_{(K)}\otimes_KK\cong X_{(K)}$. 

(iii) One has for $f:X_{(K)}\longrightarrow Y_{(K)}$
\ga{b04}{f^K:
X\otimes_kK\xrightarrow{\bar f\otimes \id} Y\otimes_kK\otimes_kK\xrightarrow{\id\otimes m} Y\otimes_kK.}
\end{lem}

Identifying $K$ with $K\otimes_kI$, we can consider $K$ as an algebra
in $\Ind\sT$ and consider the category $\Ind\sT_K$ of $K$-modules. Then
the functor ${-}^K$ factors through a functor to $\Ind\sT_K$, denoted by
the same symbol.

\begin{defn}\label{def-b03}\rm We say that
$\sT$ is flat over $K$ if the tensor product
over $K$ is exact.\end{defn}

Consequently we have the following theorem, which is analogous to Theorem \ref{thm-q09}.
\begin{thm}\label{thm-b04}Assume that the base change $\sT_{(K)}$ exists. 
For an object $U\in\sT_{(K)}$, there exists a map $U^K\otimes_kK\longrightarrow U^K$ making
$U^K$ a $K$-module. Consequently $-^K$ factors through a $K$-linear  functor to 
$\Ind\sT_K$. If $\sT$ is flat over  $K$ then this functor is a $K$-linear exact tensor functor.
\end{thm}
\begin{cor}\label{cor-b05}
Assume that $\sT$ is flat over $K$ and that 
$\sT_{(K)}$ exists. Then $\sT_{(K)}$ is equivalent to the category of triples
$(X,Y,f:X\longrightarrow Y\otimes_kK)$, $X,Y\in\sT$, $f\in\text{Mor }(\Ind\sT)$,
 and morphism defined as in \ref{sect-q10}. Consequently, any $k$-linear tensor functor
 $\omega$ from $\sT$ to a $K$-linear tensor category $\sC$ factors through
 a $K$-linear tensor $\omega_K$ from $\sT_{(K)}$ to $\sC$.
\end{cor}
\proof The first claim is proved analogously as in the proof of Proposition \ref{prop-q11}.
The second claim is a consequence of the first one. In deed,
 the functor $\omega$ induces a $K$-linear map
$$\omega\otimes K:\Hom_{\sT}(X,Y)\otimes_kK\longrightarrow \Hom_{\sC}(\omega(X),\omega(Y))$$
Thus, considering $f$ as an element of 
$\Hom_\sT(X,Y)\otimes_kK$ by means of \eqref{b01},
 we can define $\omega_K$ on an object $(X,Y,f:X\longrightarrow Y\otimes K)$ as the image of
$(\omega\otimes K) f:\omega(X)\longrightarrow \omega(Y)$. $\Box$

\begin{cor}\label{cor-b06}
Let $\sT$ be a Tannaka category over $k$. Then for any field $K\supset k$, the
base change $\sT_{(K)}$ exists. More precisely, let $K_1\supset K$ be such that 
there exists a fiber functor to $\vect_{K_1}$,  and let 
$G^{K^1}_k$ be the corresponding Tannaka groupoid.
Then $\sT_{K}\cong \Rep_f(G^{K^1}_{K})$.  
\end{cor}
\proof This follows from Lemma  \ref{lem-gr14} for 
and Theorem \ref{thm-b04}.
$\Box$

\subsection{The existence of base change by finite separable field extensions}\label{sect-b08}
Let $k\subset K$ be a finite separable field extension. Then we can repeat the argument of the proof of
Theorem \ref{thm-q18} to obtain the following result.
\begin{thm} For any finite separable field extension $K$ of $k$, the base change $\sT_{(K)}$ exists.
More precisely, $\sT_{(K)}$ can be defined as the category $\sT_K$ of $K$-modules in $\sT$.
\end{thm}
\section{Quotient category by an arbitrary Tannaka subcategory}\label{sect-g}

\subsection{}\label{sect-g07} Let $K\supset k$ be a field extension and
let $(\sQ,\fq)$ be a $K$-quotient of $\sT$ with $\sS$ being the
invariant category and assume that the base change $\sT_{(K)}$ of $\sT$
exists and that $\sT$ is flat over $K$. Then $\fq$ induces a non-neutral fiber
functor $\omega$ for $\sS$. 
By Corollary \ref{cor-b05}, $\fq$ factors through a $K$-linear 
tensor functor $\fq_K:\sT_{(K)}\longrightarrow \sQ$, which can easily be shown to be a
 quotient functor of tensor categories
over $K$ and the invariant category of $\fq_K$ is equivalent to $\sS_K$.

Let $\sO$ denote the function algebra of the Tannaka groupoid $A^K_k$ of $(\sS,\omega).$
Consider $\sO$ as an object of $\Ind\sS$ by means of the right regular coaction
$\Delta:\sO\longrightarrow \sO{}_s\otimes_t\sO$. The image of $K=I\otimes_kK$ under
$\omega$ is $K\otimes K$. Thus $\sO$ is a $K$-module in the sense of
section \ref{sect-b}, through the functor $\omega$ this means that
$\sO$ is a $K\otimes K$ modules by means of the 
map $t\otimes s$. Further $\sO$ is a $K$-algebra,
i.e. an algebra in $\Ind\sT_{(K)}$,
which means that $\sO$ is an algebra over $K\otimes K$. We define by $\Mod_{\sO,K}$
the category of $\sO$-modules in $\Ind\sT_{(K)}$, whose objects are thus
$K$-object equipped with an action of $\sO$. This is a tensor category with
respect to the tensor product over $\sO$.
\begin{defn}\label{def-g08}\rm  Let $(\sS,\omega:\sS\longrightarrow \vect_K)$  be a Tannaka subcategory of $\sT$
and $\sO$ be the corresponding function algebra considered as object of $\Ind\sT$.
We say that $\sT$ is flat over $\sS$ if it is flat over $K$ and
in $\Ind\sT_{(K)}$ the tensor product over $\sO$ is exact.
\end{defn}
\begin{thm}\label{thm-g09}Let $\sT$ be a rigid tensor category over $k$.
Let $(\sS,\omega:\sS\longrightarrow\vect_K)$ be a Tannaka subcategory, which is closed
under taking subquotients. Assume that 
\begin{itemize}\item[(i)] $\sT$ is flat over $K$ and the base change $\sT_{(K)}$ exists;
\item[(ii)] $\sT$ is flat over $\sS$ and  a $K$-quotient
$(\sQ,\fq)$ of $\sT$ by $\sS$ exists. \end{itemize}
Then the right adjoint functor $\fp$ to $\fq$ induces an
exact tensor functor to the category $\Mod_{\sO,K}$.
\end{thm}
\proof Denote by $\sO_K$ the function algebra of $A^K_K$ - the diagonal subgroup
of $A^K_k$, where $A^K_k$ is the Tannaka groupoid of $\sS$. Let $\fp_K$ denote the
right adjoint functor to $\fq_K:\Ind\sT_{(K)}\longrightarrow \Ind\sQ$. Then according to \ref{thm-q09}, $\fp_K$ is 
an exact tensor functor to the category of $\sO_K$-modules in $\Ind\sT_{(K)}$.

On the other hand, since $\sS\cong\Rep_f(A^K_k)$, 
$\sS_{(K)}\cong\Rep(A^K_K)$, the functor
$-^K$ is in fact the functor $\ind_{\Delta_K}$ in the situation of \ref{sect-gr13}. Thus 
$$(\sO_K)^K\cong \sO_K\Box_{\sO_K}\sO\cong\sO.$$
This isomorphism gives us a morphism
$$(U\otimes_{\sO_K} V)^K\longrightarrow U^K\otimes_{\sO}V^K.$$
We don't know if this morphism is an isomorphism for an arbitrary pair $U,V$ of $\sO_K$-modules in $\sT_{(K)}$
but we know it is for $U,V$ in the image of $\fq_K$, thanks to the assumption that
$\sT$ is flat over $\sS$. In deed, for $X\in\sT$, we have $\fp\fq(X)\cong X\otimes\sO$, using
the method as in the proof of \ref{thm-q09} we deduce the required isomorphism. $\Box$.

\subsection{Quotient by an \'etale finite subcategory} \label{sect-g04}
We say that a Tannaka category $\sS$ is
\'etale finite if its Tannaka groupoid is finite and reduced (with respect to the structure map
$(s,t)$). In this case there exists a fiber functor for $\sS$
with value in $\vect_K$, where $K$ is  a finite separable extension of $k$.
According to \ref{sect-b08}, $\sT$ is flat over $K$ and the base change $\sT_{(K)}$ of $\sT$ exists
and is equivalent to the category of $K$-modules in $\sT$. Further, since the category
$\sS_{(K)}$ is finite, by \ref{thm-q18} the quotient of $\sT{(K)}$ by $\sS_{(K)}$ exists and can be given as
the category of $\sO_K$-modules in $\sT{(K)}$. This is the $K$-quotient
of $\sT$ by $\sS$ we look for.

 \subsection{Application to representation theory of groupoids}
\label{sect-g10}Let $f:G^K_k\longrightarrow A^{K_0}_k$ be a homomorphism
of transitive groupoids and $L^K_{K_0}$ be the kernel of $f$. Then $\Rep_f(L^K_{K_0})$
is equivalent to the category whose objects are triples $(U,V,f:U\longrightarrow V\otimes\sO(A))$,
where $U,V\in\Rep_f(G)$, $f$ is $G$-equivariant ($G$ acts diagonally on $V\otimes\sO(A)$).

Alternatively, since morphisms $f: X\longrightarrow Y\otimes\sO$ are in 1-1 correspondence
with morphism $f^\sharp:I\longrightarrow X^\vee\otimes Y\otimes\sO$. As noticed in \ref{lem-q01}, such a
morphism $f^\sharp$ corresponds to an element of $\omega((X^\vee\otimes Y)_\sS)$.
Thus objects of $\sQ$ can be characterized as triples
$(X,Y,f\in\omega((X^\vee\otimes Y )_\sS))$.

Another consequence of the above theorem which may be useful in checking whether
a sequence of groupoids $L^K_{K_0}\longrightarrow G^K_k\longrightarrow A^{K_0}_k$ is exact is the following
\begin{cor}\label{cor-g11} Assume we are given a sequence of homomorphisms of groupoids
$L^K_{K_0}\xrightarrow q G^K_k\xrightarrow f A^{K_0}_k$ with $fq$ trivial. Then the sequence
is exact in the sense that $q$ is the kernel of $f$ iff the representation categories of these
groups satisfy the condition (i) and (ii) of Definition \ref{def-def02}
\end{cor}
\subsection{Quenstions}\label{rm-g12}The following questions are open to us:
\begin{enumerate}\item Is $\sT$ always flat over $\sS$, $K$?
\item Assume that $\sT$ is flat over $\sS$ (and $K$), does a quotient exist, in particular,
does the base change to $K$ exist?
\item Assume that $\sT$ is flat over $\sS$ (and $K$), does the functor $\fp$ induce
an equivalence of tensor categories between $\Ind\sQ$ and $\Mod_{\sO,K}$?
\end{enumerate}

Finally we mention
a related question raised by Deligne, namely whether the the fundamental group
of $\sT$ is flat over $\sS$, see \cite{DeGroth} for definition of the fundamental
group of a rigid tensor category.

\end{document}